\documentclass[10pt,twoside]{article}
\usepackage{Latex-document}
\usepackage{psfig}
\usepackage{amssymb}

\title{\bf Complex Hyperbolic Triangle Groups\thanks{Supported by
N.S.F. Research Grant DMS-0072706.}}
\author{Richard Evan Schwartz\vspace*{-0.5cm}\thanks{Department of Mathematics,
University of Maryland, College park, MD 20742, USA. E-mail:
res@math.umd.edu}}
\date{\vspace{-8mm}}

\newtheorem{theorem}{Theorem}[section]

\newtheorem{conjecture}[theorem]{Conjecture}

\def\C{\mbox{\boldmath{$C$}}}%
\def\H{\mbox{\boldmath{$H$}}}%
\def\P{\mbox{\boldmath{$P$}}}%
\def\R{\mbox{\boldmath{$R$}}}%
\def\Z{\mbox{\boldmath{$Z$}}}%

\begin{document}

\maketitle

\thispagestyle{first} \setcounter{page}{339}

\begin{abstract}

\vskip 3mm

The theory of complex hyperbolic discrete groups is still in its childhood but promises to grow into a rich
subfield of geometry.  In this paper I will discuss some recent progress that has been made on complex hyperbolic
deformations of the modular group and, more generally, triangle groups. These are some of the simplest nontrivial
complex hyperbolic discrete groups.  In particular, I will talk about my recent discovery of a closed real
hyperbolic 3-manifold which appears as the manifold at infinity for a complex hyperbolic discrete group.

\vskip 4.5mm


\end{abstract}

\vskip 12mm

\section{Introduction}

\vskip-5mm \hspace{5mm}

A basic problem in geometry is the {\it deformation problem\/}.
One starts with a finitely generated group $\Gamma$, a Lie group
$G_1$, and a larger Lie group $G_2 \supset G_1$. Given a discrete
embedding $\rho_0: \Gamma \to G_1$ one asks if $\rho_0$ fits
inside a family $\rho_t: \Gamma \to G_2$ of discrete embeddings.
Here {\it discrete embedding\/} means an injective homomorphism
onto a discrete set.

A nice setting for the deformation problem is the case when $G_1$
and $G_2$ are isometry groups of rank one symmetric spaces, $X_1$
and $X_2$, and $\Gamma$ is isomorphic to a lattice in $G_1$. If
$X_1=\H^2$, the hyperbolic plane, and $X_2=\H^3$, hyperbolic
$3$-space, then we are dealing with the classic and well-developed
theory of quasifuchsian groups.

The $(p,q,r)$-{\it reflection triangle group\/} is possibly the
simplest kind of lattice in Isom$(\H^2)$.  This group is generated
by reflections in the sides of a geodesic triangle having angles
$\pi/p$, $\pi/q$, $\pi/r$ (subject to the inequality
$1/p+1/q+1/r<1$.)  We allow the possibility that some of the
integers are infinite. For instance, the $(2,3,\infty)$-reflection
triangle group is commensurable to the classical modular group.

The reflection triangle groups are rigid in Isom$(\H^3)$, in the
sense that any two discrete embeddings of the same group are
conjugate. We are going to replace $\H^3$ by $\C\H^2$, the complex
hyperbolic plane. In this case, we get nontrivial deformations.
These deformations provide an attractive problem, because they
furnish some of the simplest interesting examples in the still
mysterious subject of complex hyperbolic deformations.  While some
progress has been made in understanding these examples, there is
still a lot unknown about them.

In \S 2 we will give a rapid introduction to complex hyperbolic
geometry.  In \S 3 we will explain how to generate some complex
hyperbolic triangle groups. In \S 4 we will survey some results
about these groups and in \S 5 we will present a more complete
conjectural picture. In \S 6 we will indicate some of the
techniques we used in proving our results.

\section{The complex hyperbolic plane}

\vskip-5mm \hspace{5mm}

The book [{8\/}] is an excellent general reference for complex
hyperbolic geometry. Here are some of the basics.

$\C^{2,1}$ is a copy of the vector space $\C^{3}$ equipped with
the Hermitian form
\begin{equation}
\langle  U, V\rangle=-u_{3} \overline v_{3}+ \sum_{j=1}^n u_j
\overline v_j.
\end{equation}
Here $U=(u_1,u_2,u_3)$ and $V=(v_1,v_2,v_3)$. A vector $V$ is
called {\it negative\/}, {\it null\/}, or {\it positive\/}
depending (in the obvious way) on the sign of $\langle V,V
\rangle$.  We denote the set of negative, null, and positive
vectors, by $N_-$, $N_0$ and $N_+$ respectively.

$\C^2$ includes in complex projective space $\C\P^2$ as the affine
patch of vectors with nonzero last coordinate. Let $[\ ]:
\C^{2,1}-\{0\} \to \C\P^2$ be the projectivization whose formula,
expressed in the affine patch, is
\begin{equation}
\label{theta} [(v_1,v_2,v_3)]=(v_1/v_3,v_2/v_3).
\end{equation}
The {\it complex hyperbolic plane\/}, $\C\H^2$, is the projective
image of the set of negative vectors in $\C^{2,1}$.  That is,
$\C\H^2=[N_-]$.  The ideal boundary of $\C\H^2$ is the unit sphere
$S^{3}=[N_0]$. If $[X],[Y] \in \C\H^n$ the complex hyperbolic
distance $\varrho([X],[Y])$ satisfies
\begin{equation}
\label{hypdist}\label{hyp} \varrho([X],[Y])=2
\cosh^{-1}\sqrt{\delta( X, Y)}; \hskip 20 pt \delta( X, Y) =\frac{
\langle X, Y\rangle \langle Y, X \rangle} {\langle X, X\rangle
\langle Y, Y\rangle}.
\end{equation}
Here $X$ and $Y$ are arbitrary lifts of $[X]$ and $[Y]$. See
[{8\/}, 77]. The distance we defined is induced by an invariant
Riemannian metric of sectional curvature pinched between $-1$ and
$-4$.  This Riemannian metric is the real part of a K\"ahler
metric.

$SU(2,1)$ is the Lie group of $\langle \ ,\  \rangle$ preserving
complex linear transformations. $PU(2,1)$ is the projectivization
of $SU(2,1)$ and acts isometrically on $\C\H^2$. The map $SU(2,1)
\to PU(2,1)$ is a $3$-to-$1$ Lie group homomorphism. The group of
holomorphic isometries of $\C\H^2$ is exactly $PU(2,1)$. The full
group of isometries of $\C\H^2$ is generated by $PU(2,1)$ and by
the antiholomorphic map $(z_1,z_2,z_3) \to (\overline
z_1,\overline z_2,\overline z_3)$.

An element of $PU(2,1)$ is called {\it elliptic\/} if it has a
fixed point in $\C\H^2$.  It is called {\it hyperbolic\/} (or {\it
loxodromic\/}) if there is some $\epsilon>0$ such that every point
in $\C\H^2$ is moved at least $\epsilon$ by the isometry. An
element which is neither elliptic nor hyperbolic is called {\it
parabolic\/}.

$\C\H^2$ has two different kinds of totally geodesic subspaces,
{\it real slices\/} and {\it complex slices\/}.  Every real slice
is isometric to $\C\H^2 \cap \R^2$ and every complex slice is
isometric to $\C\H^2 \cap \C^1$. The ideal boundaries of real and
complex slices are called, respectively, $\R$-circles and
$\C$-circles. The complex slices naturally implement the
Poincar\'e model of the hyperbolic plane and the real slices
naturally model the Klein model.  It is a beautiful feature of the
complex hyperbolic plane that it contains both models of the
hyperbolic plane.

\section{Reflection triangle groups}

\vskip-5mm \hspace{5mm}

There are two kinds of reflections in Isom$(\C\H^2)$. A {\it real
reflection\/} is an anti-holomorphic isometry conjugate to the map
$(z,w) \to (\overline z, \overline w)$. The fixed point set of a
real reflection is a real slice.   We shall not have much to say
about the explicit computation of real reflections, but rather
will concentrate on the complex reflections.

A {\it complex reflection\/} is a holomorphic isometry conjugate
to the involution $(z,w) \to (z,-w)$.  The fixed point set of a
complex reflection is a complex slice. There is a simple formula
for the general complex reflection: Let $ C \in  N_+$. Given any $
U \in \C^{2,1}$ define
\begin{equation}
\label{refl} I_{C}( U)=- U+ \frac{2 \langle  U, C\rangle}{\langle
C,
 C\rangle}   C.
\label{reflect}
\end{equation}
$I_{C}$ is a complex reflection.

We also have the formula
\begin{equation}
\label{Hcross} \label{hcross}
 U \boxtimes  V = (\overline{u_3v_2-u_2v_3},
                \ \overline{u_1v_3-u_3v_1},
                \ \overline{u_1v_2-u_2v_1}).
\end{equation}
This vector is such that $\langle  U, U
 \boxtimes  V\rangle
=\langle  V, U \boxtimes  V\rangle =0$. See [{8\/}, p. 45].

Equations \ref{refl} and \ref{Hcross} can be used in tandem to
rapidly generate triangle groups defined by complex reflections.
One picks three vectors $V_1,V_2,V_2 \in N_-$. Next, we let
$C_j=V_{j-1} \boxtimes V_{j+1}$. Indices are taken mod $3$.
Finally, we let $I_j=I_{C_j}$.   The complex reflection $I_j$
fixes the complex line determined by the points $[V_{j-1}]$ and
$[V_{j+1}]$. This, the group $\langle I_1,I_2,I_3 \rangle$ is a
complex-reflection triangle group determined by the triangle with
vertices $[V_1]$, $[V_2]$, $[V_3]$.

Here is a quick dimension count for the space of
$(p,q,r)$-triangle groups generated by complex reflections.  We
can normalize so that $[V_1]=0$. The stabilizer of $0$ in
$PU(2,1)$ acts transitively on the unit tangent space at $0$.  We
can therefore normalize so that $[V_2]=(s,0)$ where $s \in (0,1)$.
Finally, the isometries $(z,w) \to (z,\exp(i \theta) w)$ stabilize
both $[V_1]$ and $[V_2]$.  Applying a suitable isometry we arrange
that $[V_3]=(t+i u,v)$ where $t,u,v \in (0,1)$.   We cannot make
any further normalizations, so the space of triangles in $\C\H^2$
mod isometry is $4$-real dimensional. Each of the three angles
$(p,q,r)$ puts $1$ real constraint on the triangle.  For instance,
the $p$-angle places the constraint that $(I_1I_2)^p$ is the
identity. Since $4-3=1$, we see heuristically that the space of
$(p,q,r)$-complex reflection triangle groups is $1$-real
dimensional.

The argument we just gave can be made rigorous, and extends to the
case when some of the integers are infinite. (In this case the
corresponding vectors are null rather than negative.) In the
$(\infty,\infty,\infty)$-case, the parameter is the {\it angular
invariant\/} $\arg (\langle V_1,V_2 \rangle \langle V_2,V_3
\rangle \langle V_3,V_1 \rangle)$. Compare [{10\/}].

This $1$-dimensionality of the deformation space makes the
$(p,q,r)$-triangle groups an especially attractive problem to
study.  Indeed, there is a completely canonical path of
deformations. The starting point for the path of deformations is
the case when the vectors have entirely real entries. (That is,
$u=0$.) In this case, the three complex reflections stabilize the
real slice $\R^2 \cap \C\H^2$.

\section{Some results}

\vskip-5mm \hspace{5mm}

To obtain a deformation of the $(p,q,r)$-reflection triangle group
we choose a slice, either real or complex, and a triple of
reflections, either real or complex, which restrict to the
reflections in the sides of a $(p,q,r)$-geodesic triangle in the
slice. {\it A priori\/} there are $4$ possibilities, given that
the slice and the reflection types can be either real or complex.
These choices lead to different outcomes.

If we start with complex reflections stabilizing a complex slice,
the group has order $2$, because the reflections will all
stabilize the same slice.

A more interesting case involving complex slices is given by:
\renewcommand{\thesection}{\arabic{section}}
\begin{theorem}{\rm [8]}
$\rho_0: \Gamma \to $Isom$(\C\H^2)$ stabilizes a complex slice and
acts on this slice with compact quotient then any nearby
representation $\rho_t$ also stabilizes a complex slice.
\end{theorem}

Goldman's theorem applies to any co-compact lattice
$\rho_0(\Gamma)$.  In the case of triangle groups, which are rigid
in $\H^2$,  it says that any nearby representation is conjugate
the original.  In contrast:

\begin{theorem} {\rm [4, 12]} There is a $1$-parameter
family $\rho_t(\Gamma(2,3,\infty))$ of discrete faithful
representations of the modular group having the property that
$\rho_0$ stabilizes a real slice and $\rho_1$ stabilizes a complex
slice.   For every parameter the generators are real reflections.
\end{theorem}

Thus, in the case of non-cocompact triangle groups, two of the
remaining $3$ cases can be {\it connected\/}. In their paper,
Falbel and Koseleff claim that their technique works for
$\Gamma(p,q,\infty)$ when $\max(p,q)=4$.  For higher values of $p$
and $q$ it is not known what happens.

The remaining case occurs when we start with complex reflections
stabilizing a real slice. This is the case we discussed in the
previous section.  Henceforth we restrict our attention to this
case.

Goldman and Parker introduced this topic and studied the case of
the ideal triangle group $\Gamma(\infty,\infty,\infty)$. They
found that there is a $1$-real parameter family of non-conjugate
representations, $\{\rho_t$, $t \in (-\infty,\infty)\}$. Once
again $\rho_0$ stabilizes a real slice. Paraphrasing their more
precise formulation:

\begin{theorem}{\rm [10]} There are symmetric
neighborhoods $I \subset J$ of $0$ such that $\rho_t$ is discrete
and faithful if $t \in I$ and not both discrete and faithful if $t
\not \in J$.
\end{theorem}

$J$ consists of the parameter values $t$ such that the element
$\rho_t(I_1I_2I_3)$ is not an elliptic element.  For $t \not \in
J$, this element is elliptic.  If it has finite order then the
representation is not faithful;  if it has infinite order then the
representation is not discrete. The (very slightly) smaller
interval $I$ is the interval for which their proof works. They
conjectured that $\rho_t$ should be discrete and faithful iff $t
\in J$.

We proved the Goldman-Parker conjecture, and sharpened it a bit.

\begin{theorem}{\rm [16]}
\label{one} $\rho_t$ is discrete and faithful if and only if $t
\in J$.  Furthermore, $\rho_t$ is indiscrete if $t \not \in J$.
\end{theorem}

The group $L=\rho_s(\Gamma(\infty,\infty,\infty))$, when $s \in
\partial J$ is especially beautiful. We call this group the {\it
last ideal triangle group\/}. (There are really two groups, one
for each endpoint of $J$, but these are conjugate.) This group
seems central in the study of complex hyperbolic deformations of
the modular group.  For instance, Falbel and Parker recently
discovered that $L$ arises as the endpoint of a certain family of
deformations of the modular group, using real reflections.  See
[{5\/}] for details.

Recall that $L$, like all discrete groups, has a {\it limit set\/}
$\Omega(L) \subset S^3$ and a {\it domain of discontinuity\/}
$\Delta(L)=S^3-\Omega(L)$.   The quotient $\Delta(L)/L$ is a
$3$-dimensional orbifold, commonly called the {\it orbifold at
infinity\/}.

\begin{theorem}{\rm [17]}
\label{two} $\Delta(L)/L$ is commensurable to the Whitehead link
complement.
\end{theorem}

The Whitehead link complement is a classic example of a finite
volume hyperbolic $3$-manifold.  The surprise in the above result
is that a real hyperbolic $3$-manifold makes its appearance in the
context of complex hyperbolic geometry.

One might wonder about analogues of Theorem \ref{one} for other
triangle groups.  Below we will conjecture that the space of
discrete embeddings is a certain interval. In his thesis [{22\/}],
Justin Wyss-Gallifent studied some special cases of this question.
He made a very interesting discovery concerning the $(4,4,\infty)$
triangle group:

\begin{theorem}{\rm [22]}
Let $S$ be the set of parameters $t$ for which the representation
$\rho_t(\Gamma(4,4,\infty))$ is discrete (but not necessarily
injective).  Then $S$ contains isolated points and, in particular,
is not an interval.
\end{theorem}

There seems to be an interval $J$ of discrete embeddings and,
outside of $J$, an extra countable sequence $\{t_j\}$ of
parameters for which $\rho_{t_j}$ is discrete but not an
embedding. This sequence accumulates on the endpoints of $J$.

Motivated by [{17\/}] I wanted to produce a discrete complex
hyperbolic group whose orbifold at infinity was a closed
hyperbolic $3$-manifold. The extra representations found by
Wyss-Gallifent seemed like a good place to start.  Unfortunately,
there is a cusp built into the representations of the
$(4,4,\infty)$ triangle groups.

Instead, I considered the $(4,4,4)$-groups, and found that the
extra discrete deformations exist. $\rho_t(\Gamma(4,4,4))$ seems
to be discrete embedding iff all the elements of the form
$\rho_t(I_iI_jI_iI_k)$ are not elliptic. Here $i,j,k$ are meant to
be distinct. (For all these parameters, the element
$\rho_t(I_iI_jI_k)$ is still a loxodromic element.) There is a
countable collection $t_5,t_6,...$ of parameters such that
$\rho_{t_j}(I_iI_jI_iI_k)$ has order $j$. All these
representations seem discrete. For ease of notation we set
$\rho_j=\rho_{t_j}$.

For $j=5,6,7,8,12$ we can show by arithmetic means that $\rho_j$
is discrete.   The representation $\rho_5$ was too complicated for
me to analyze and $\rho_6$ has a cusp. The simplest remaining
candidate is $\rho_7$.

\begin{theorem}{\rm [18]}
\label{three} $G=\rho_7(\Gamma(4,4,4))$ is a discrete group. The
orbifold at infinity $\Delta(G)/G$ is a closed hyperbolic
$3$-orbifold.
\end{theorem}
\renewcommand{\thesection}{\arabic{section}.}

In the standard terminology, $\Delta(G)/G$ is the orbifold
obtained by labelling the braid $(AB)^{15}(AB^{-2})^3$ with a $2$.
Here $A$ and $B$ are the standard generators of the $3$-strand
braid group.

A {\it spherical CR structure\/} on a $3$-manifold is a system of
coordinate charts into $S^3$ whose transition functions are
restrictions of complex projective transformations. Kamishima and
Tsuboi [{13\/}] produced examples of spherical CR structures on
Seifert fibered $3$-manifolds, but our example in theorem
\ref{three} gives the only known spherical CR structure on a
closed hyperbolic $3$-manifold. We think that Theorem \ref{three}
holds for all $j=8,9,10...$.

Concerning the specific topic of triangle groups generated by
complex reflections, I think that not much else is known. Recently
a lot of progress has been made in understanding triangle groups
generated by real reflections. See [{3\/}] and [{4\/}]. There has
been a lot of other great work done recently on complex hyperbolic
discrete groups, for instance [{1\/}], [{2\/}], [{9\/}], [{20\/}],
[{21\/}].  Also see the references in Goldman's book [f{8\/}].

\section{A conjectural picture}

\vskip-5mm \hspace{5mm}

We will consider the $1$-parameter family $\rho_t(p,q,r)$ of
representations of the $(p,q,r)$-reflection triangle group, using
complex reflections.  We arrange that $\rho_0$ stabilizes a real
slice. We choose our integers so that $p \leq q \leq r$.  We let
$I_p$, $I_q$, $I_r$ be the generators of the reflection triangle
group.  The notation is such that $I_p$ is the reflection in the
side of the triangle opposite $p$, etc. Define
\begin{equation}
W_A=I_pI_rI_qI_r; \hskip 30 pt W_B=I_pI_qI_r.
\end{equation}

\renewcommand{\thesection}{\arabic{section}}
\begin{conjecture}
\label{discrete} The set of $t$ for which $\rho_t(p,q,r)$ is a
discrete embedding is the closed interval consisting of the
parameters $t$ for which neither $\rho_t(W_A)$ nor $\rho_t(W_B)$
is elliptic.
\end{conjecture}

We call the interval of Conjecture \ref{discrete} the {\it
critical interval\/}.

We say that the triple $(p,q,r)$ has {\it type $A$\/} if the
endpoints of the critical interval correspond to the
representations when $W_A$ is a parabolic element.  In other
words, $W_A$ becomes elliptic before $W_B$.  We say otherwise that
$(p,q,r)$ has {\it type $B$\/}.

\begin{conjecture}
The triple $(p,q,r)$ has type A if $p<10$ and type B if $p>13$.
\end{conjecture}

The situation is rather complicated when $p \in \{10,11,12,13\}$.
Our Java applet [{19\/}] lets the user probe these cases by hand,
though the roundoff error makes a few cases ambiguous.  The extra
deformation, which was the subject of Theorem \ref{three}, seems
part of a more general pattern.

\begin{conjecture}
\label{extra} If $(p,q,r)$ has type A then there is a countable
collection of parameters $t_1,t_2,t_3...$ for which
$\rho_{t_j}(p,q,r)$ is infinite and discrete but not injective. If
$(p,q,r)$ has type B then all infinite discrete representations
$\rho_t(p,q,r)$ are embeddings and covered by Conjecture
\ref{discrete}.
\end{conjecture}

\noindent The {\it proviso\/} about the infinite image arises
because there always exists an extremely degenerate representation
of $\Gamma(p,q,r)$ onto $\Z/2$.  The generators are all mapped to
the same complex reflection.

In summary, there seems to be a critical interval $I$, such the
representations $\rho_t(p,q,r)$ are discrete embeddings iff $t \in
I$.  Depending on the endpoints of $I$, there are either no
additional discrete representations, or a countable collection of
extra discrete representations.

It is interesting to see what happens as $t$ moves to the boundary
of $I$ from within $I$.  We observed a certain kind of
monotonicity to the way the representation varies. Let $\Gamma$ be
the abstract $(p,q,r)$ triangle group. For any word $W \in
\Gamma$, let $W_t=\rho_t(W)$.   We will concentrate on the case
when $W$ is an infinite word. For $t \in I$, the element $W_t$ is
(conjecturally) either a parabolic or loxodromic. Let
$\lambda(W_t)$ be the translation length of $W_t$.

\begin{conjecture}
\label{mono} As $t$ increases monotonically from $0$ to $\partial
I$, the quantity $\lambda(W_t)$ decreases monotonically for all
infinite words $W$.
\end{conjecture}
\renewcommand{\thesection}{\arabic{section}.}

Conjecture \ref{mono} is closely related to some conjectures of
Hanna Sandler [{15\/}] about the behavior of the trace function in
the ideal triangle case. I think that there is some fascinating
algebra hiding behind the triangle groups$-$in the form of the
behavior of the trace function---but so far it is unreachable.

\section{Some techniques of proof}

\vskip-5mm \hspace{5mm}

If $G \subset $Isom$(X)$, one can try to show that $G$ is discrete
by constructing a {\it fundamental domain\/} for $G$.   One looks
for a set $F \subset X$ such that the orbit $G(F)$ {\it tiles\/}
$X$.  This means that the translates of $F$ only intersect $F$ in
its boundary. The Poincar\'e theorem [{B\/}, \S 9.6] gives a
general method for establishing the tiling property of $F$ based
on how certain elements of $G$ act on $\partial F$.

When $X=\H^n$, one typically builds fundamental domains out of
polyhedra bounded by totally geodesic codimension-$1$ faces. When
$X=\C\H^n$, the situation is complicated by the absence of totally
geodesic codimension-$1$ subspaces.  The most natural replacement
is the {\it bisector\/}.  A bisector is the set of points in
$\C\H^n$ equidistant between two given points. Mostow [{14\/}]
used bisectors in his analysis of some exceptional non-arithmetic
lattices in Isom$(\C\H^2)$, and Goldman studied them extensively
in [{8\/}]. (See Goldman's book for additional references on
papers which use bisectors to construct fundamental domains.)

My point of view is that there does not seem to be a ``best'' kind
surface to use in constructing fundamental domains in complex
hyperbolic space. Rather, I think that one should be ready to
fabricate new kinds of surfaces to fit the problem at hand.  It
seems that computer experimentation often reveals a good choice of
surface to use.  In what follows I will give a quick tour of
constructive techniques.

Consider first the deformations $G_t=\rho_t(\infty,\infty,\infty)$
of the ideal triangle group, introduced in [{10\/}]. According to
[{16\/}] these groups are discrete for $t \in [0,\tau]$. Here
$\tau$ is the {\it critical parameter\/} where the product of the
generators is parabolic. It is convenient to introduce the {\it
Clifford torus\/}.  Thinking of $\C\H^2$ as the open unit ball in
$\C^2$, the Clifford torus is the subset $T=\{|z|=|w|\} \subset
S^3$.  Amazingly $T$ has $3$ foliations by $\C$-circles: The {\it
horizontal foliation\/} consists of $\C$-circles of the form
$\{(z,w)|z=z_0\}$. The {\it vertical foliation\/} consists of
$\C$-circles of the form $\{(z,w)|w=w_0\}$. The {\it diagonal
foliation\/} consists of $\C$-circles having the form $\{(z,w)|\
z=\lambda_0 w\}$.

Recall that $G_t$ is generated by $3$ complex reflections.  Each
of these reflections fixes a complex slice and hence the bounding
$\C$-circle. One can normalize so that the three fixed
$\C$-circles lie on the Clifford torus, one in each of the
foliations.  Passing to an index $2$ subgroup, we can consider a
group generated by $4$ complex reflections: Two of these
reflections, $H_1$ and $H_2$, fix horizontal $\C$-circles $h_1$
and $h_2$  and the other two, $V_1$ and $V_2$, fix vertical
$\C$-circles $v_1$ and $v_2$.

The ideal boundary of a bisector is called a {\it  spinal
sphere.\/} This is an embedded $2$-sphere which is foliated by
$\C$-circles (and also by $\R$-circles.)  We can find a
configuration of $4$-spinal spheres $S(1,v)$, $S(2,v)$, $S(1,h)$
and $S(2,h)$. Here $S(j,v)$ contains $v_j$ as part of its
foliation and $S(j,h)$ contains $h_j$ as part of its foliation.
The map $H_j$ stabilizes $S(j,h)$ and interchanges the two
components of $S^3-S(j,h)$.  Analogous statements apply to the
$V$s.

The two spheres $S(h,j)$ are contained in the closure of one
component of $S^3-T$ and the two spheres $S(v,j)$ are contained in
the closure of the other.  When the parameter $t$ is close to $0$
these spinal spheres are all disjoint from each other, excepting
tangencies, and form a kind of necklace of spheres.  Given the way
the elements $H_j$ and $V_j$ act on our necklace of spheres, we
see that we are dealing with the usual picture associated to a
{\it Schottky group\/}.  In this case the discreteness of the
group is obvious.

As the parameter increases, the two spinal spheres $S(v,1)$ and
$S(v,2)$ collide. Likewise, $S(h,1)$ and $S(h,2)$ collide.
Unfortunately, the collision parameter occurs before the critical
parameter. For parameters larger than this collision parameter, we
throw out the spinal spheres and look at the action of $G$ on the
Clifford torus itself. (This is not the point of view taken in
[{10\/}] but it is equivalent to what they did.)

Let $H$ be the subgroup generated by the reflections $H_1$ and
$H_2$.  One finds that the orbit $H(T)$ consists of translates of
$T$ which are disjoint from each other except for forced
tangencies.  Even though $H$ is an infinite group, most of the
elements in $H$ move $T$ well off itself, and one only needs to
take care in checking a short finite list of words in $H$.   Once
we know how $H$ acts on $T$ we invoke a variant of the {\it
ping-pong lemma\/} to get the discreteness.

At some new collision parameter, the translates of the Clifford
torus collide with each other. Again, the collision parameter
occurs before the critical parameter.  This is where the work in
[{16\/}] comes in. I define a new kind of surface called a {\it
hybrid cone\/}. A hybrid cone is a certain surface foliated by
arcs of $\R$-circles.  These arcs make the pattern of a fan: Each
arc has one endpoint on the arc of a $\C$-circle and the other
endpoint at a single point common to all the arcs. I cut out two
triangular patches on the Clifford torus and replace each patch by
a union of three hybrid cones. Each triangular patch is bounded by
three arcs of $\C$-circles; so that the hybrid cones are formed by
connecting these exposed arcs to auxilliary points using arcs of
$\R$-circles. In short, I put some dents into the Clifford torus
to make it fit better with its $H$-translates, and the I apply the
ping-pong lemma to the dented torus.

I also use hybrid cones in [{17\/}], to construct a natural
fundamental domain in the domain of discontinuity $\Delta(L)$ for
the last ideal triangle group $L$.  In this case, the surfaces fit
together to make three topological spheres, each tangent to the
other two along arcs of $\R$-circles. The existence of this
fundamental domain lets me compute explicitly that $\Delta(L)/L$
is commensurable to the Whitehead link complement.

Falbel and Zocca [{6\/}] introduce related surfaces called
$\C$-spheres, which are foliated by $\C$-circles. These surfaces
seem especially well adapted to groups generated by real
reflections.  See [{3\/}] and [{4\/}].  Indeed, Falbel and Parker
construct a different fundamental domain for $L$ using
$\C$-spheres.  See [{5\/}].

To prove Theorem \ref{three} in [{18\/}] I introduce another
method of constructing fundamental domains. My proof revolves
around the construction of a simplicial complex $Z\subset
\C^{2,1}$.  The vertices of $Z$ are canonical lifts to $\C^{2,1}$
of fixed points of certain elements of the group
$G=\rho_7(\Gamma(4,4,4))$. The tetrahedra of $Z$ are Euclidean
convex hulls of various $4$-element subsets of the vertices.
Comprised of infinitely many tetrahedra, $Z$ is invariant under
the element $I_2I_1I_3$.  Modulo this element $Z$ has only
finitely many tetrahedra.

Recall that $[\ ]$ is the projectivization map. Let $[Z_0]=[Z]
\cap S^3$. I deduce the topology of the orbifold at infinity by
studying the topology of $[Z_0]$. To show that my analysis of the
topology at infinity is correct, I show that one component $F$ of
$\C\H^2-[Z]$ has the {\it tiling property\/}: The $G$-orbit of $F$
tiles $\C\H^2$. Now, $Z$ is an essentially combinatorial object,
and it not too hard to analyze the combinatorics and topology of
$Z$ in the abstract.  The hard part is showing that the map $Z \to
[Z]$ is an embedding.  Assuming the embedding, the combinatorics
and topology of $Z$ are reproduced faithfully in $[Z]$, and I
invoke a variant of the Poincar\'e theorem.

After making some easy estimates, my main task boils down to
showing that the projectivization map $[\ ]$ is injective on all
pairs of tetrahedra within a large but finite portion of $Z$.
Roughly, I need to check about $1.3$ million tetrahedra. The sheer
number of checks forces us to bring in the computer. I develop a
technique for proving, with rigorous machine-aided computation,
that $[\ ]$ is injective on a given pair of tetrahedra.

A novel feature of my work is the use of computer experimentation
and computer-aided proofs.  This feature is also a drawback,
because it only allows for the analysis of examples one at a time.
To make this analysis automatic I would like to see a kind of
marriage of complex hyperbolic geometry and computation.  On the
other hand, I would greatly prefer to see some theoretical
advances in discreteness-proving which would eliminate the
computer entirely.

\label{lastpage}


\begin{thebibliography}{aa}
\bibitem{041101} D. Allcock, J. Carlson, D. Toledo, {\it The Moduli
Space of Cubic Threefolds\/}, J. Alg. Geom. (to appear).
\bibitem{041102} P. Deligne and G.D. Mostow, {\it Commensurabilities
among Lattices in $PU(1,n)$\/}, Annals of Mathematics Studies {\bf
132\/}, Princeton University Press (1993).
\bibitem{041103} E. Falbel and P.-V. Koseleff, {\it Flexibility of
the Ideal Triangle Group in Complex Hyperbolic Geometry\/},
Topology {\bf 39(6)\/} (2000), 1209--1223.
\bibitem{041104} E. Falbel and P.-V. Koseleff, {\it A Circle of
Modular Groups\/}, preprint 2001.
\bibitem{041105} E. Falbel and J. Parker, {\it The Moduli Space of
the Modular Group in Complex Hyperbolic Geometry\/}, Math.
Research Letters (to appear).
\bibitem{041106} E. Falbel and V. Zocca, {\it A Poincar\'e's
Fundamental Polyhedron Theorem for Complex Hyperbolic
Manifolds\/}, J. reine angew Math. {\bf 516\/} (1999), 133--158.
\bibitem{041107} W. Goldman {\it Representations of fundamental groups
of surfaces\/}, in ``Geometry and Topology, Proceedings,
University of Maryland 1983--1984'', J.\ Alexander and J.\ Harer
(eds.), Lecture Notes in Math. Vol.~1167 (1985), 95--117.
\bibitem{041108} W. Goldman, {\it Complex Hyperbolic Geometry\/},
Oxford Mathematical Monographs, Oxford University Press, (1999).
\bibitem{041109} W. Goldman, M. Kapovich and B. Leeb, {\it Complex
Hyperbolic Surfaces Homotopy Equivalent to a Riemann surface\/},
Communications in Analysis and Geometry {\bf 9\/} (2001), 61--95.
\bibitem{041110} W. Goldman and J. Parker, {\it Complex Hyperbolic
Ideal Triangle Groups\/}, J. reine agnew Math. {\bf 425\/} (1992),
71--86.
\bibitem{041111} W. Goldman and J. Millson, {\it Local Rigidity of
Discrete Groups Acting on Complex Hyperbolic Space\/}, Inventiones
Mathematicae {\bf 88\/} (1987), 495--520.
\bibitem{041112} N. Gusevskii and J.R. Parker, {\it Complex Hyperbolic
Representations of Surface Groups and Toledo's Invariant\/},
preprint (2001).
\bibitem{041113} Y. Kamishima and T. Tsuboi, {\it CR Structures on
Seifert Manifolds\/}, Invent. Math. {\bf 104\/} (1991) 149--163.
\bibitem{041114} G.D. Mostow, {\it On a Remarkable Class of Polyhedra
in Complex Hyperbolic Space\/}, Pac. Journal of Math {\bf 86\/}
(1980) 171--276.
\bibitem{041115} H. Sandler, {\it Trace Equivalence in $SU(2,1)$\/},
Geo Dedicata {\bf 69\/} (1998) 317--327.
\bibitem{041116} R. E. Schwartz, {\it Ideal Triangle Groups, Dented Tori, and
Numerical Analysis\/}, Annals of Math {\bf 153\/} (2001).
\bibitem{041117} R.E. Schwartz, {\it Degenerating the Complex
Hyperbolic Ideal Triangle Groups\/}, Acta Mathematica {\bf 186\/}
(2001).
\bibitem{041118} R. E. Schwartz, {\it Real Hyperbolic on the Outside,
Complex Hyperbolic on the Inside\/}, Invent. Math (to apear).
\bibitem{041119} R. E. Schwartz, {\it Applet 29\/} (2001)
http://www.math.umd.edu/$\widetilde{\hskip 10 pt}$res.
\bibitem{041120} Y. Shalom, {\it Rigidity, Unitary Representations of
Semisimple Groups, and Fundamental Groups of Manifolds with Rank
One Transformation Group\/}, Annals of Math {\bf 152\/} (2000)
113--182.
\bibitem{041121} D. Toledo, {\it Representations of Surface Groups on
Complex Hyperbolic Space\/}, Journal of Differential Geometry {\bf
29\/} (1989) 125--133.
\bibitem{041122} J. Wyss-Gallifent, {\it Discreteness and
Indiscreteness Results for Complex Hyperbolic Triangle Groups\/},
Ph.D. Thesis, University of Maryland (2000).
\end{thebibliography}
\end{document}